\documentclass[12pt]{amsart}
\usepackage{amsmath,amssymb}

\begin{document}

\newtheorem{thm}{Theorem}%[section]
\newtheorem{lem}{Lemma}[section]
\newtheorem{prop}[thm]{Proposition}%[section]
\newtheorem{cor}[thm]{Corollary}%[section]
%\newtheorem{conj}[thm]{Conjecture}
%\newtheorem{def}[thm]{Definition}%[section]

%%%%%%%
\theoremstyle{definition}
\newtheorem{defn}{Definition}%[section]
\newtheorem{conj}{Conjecture}

%%%%%%%
\theoremstyle{remark}
\newtheorem{remark}{Remark}%[section] % \renewcommand{\theremark}{}
\newtheorem{note}{Note}[section]
%%%%%%%%%

\newcommand{\Z}{{\mathbb Z}}
\newcommand{\Q}{{\mathbb Q}}
\newcommand{\R}{{\mathbb R}}
\newcommand{\C}{{\mathbb C}}
\newcommand{\supp}{\operatorname{supp}}
\newcommand{\re}{\operatorname{Re}}
\newcommand{\im}{\operatorname{Im}}
\newcommand{\dist}{\operatorname{dist}}
\newcommand{\vol}{\operatorname{vol}}
\newcommand{\ord}{\operatorname{ord}}

\newcommand{\HH}{\mathbb H} % upper halfplane
\newcommand{\BG}{B_k^{(\Gamma)}} % Bergmann kernel for Gamma

\newcommand{\Wr}{\operatorname{Wr}} % Wronskian

\newcommand{\SH}{ {\mathcal S\mathcal H} } % subharmonic
\newcommand{\SHO}{{ {\mathcal S\mathcal H} }(\Omega)}

\newcommand{\toweak}{\xrightarrow[]{w}} % weak convergence

\newcommand{\conv}{\ast} % convolution

\newcommand{\dV}{dV} % Poincare volume 

\newcommand{\ZF}{Z} % zeros measure

\newcommand{\laplaceh}{\mathcal L} % hyperbolic laplacian

\title[distribution of zeros of modular forms]
{On the asymptotic distribution of zeros of modular forms}

\author{Ze\'ev Rudnick}

\address{Raymond and Beverly Sackler School of Mathematical Sciences,
Tel Aviv University, Tel Aviv 69978, Israel 
({\tt rudnick@math.tau.ac.il})}

\date{June 14, 2005}

\thanks{Partially supported by grant number 2002088  
 from the United States-Israel Binational Science
 Foundation (BSF), Jerusalem, Israel}

 \begin{abstract}
  We study the distribution of zeros of holomorphic modular
  forms. Assuming the Generalized Riemann Hypothesis we show that the
  zeros of Hecke eigenforms for the modular group become
  equidistributed with respect to the hyperbolic  measure on the
  modular domain as the weight grows.
 \end{abstract}

\maketitle

%\input{newintro}

% newintro 31.5.2005
\section{Introduction}

\subsection{}
 Our purpose in this note is to
study the limiting distribution of zeros of modular forms.  
We review some definitions:
A modular form of weight $k$ for $SL_2(\Z)$
is a holomorphic function   on the upper half-plane $\HH$,
transforming as $f(\frac{az+b}{cz+d}) = (cz+d)^k f(z)$, for all $
(\begin{smallmatrix}  a&b\\c&d \end{smallmatrix})\in SL_2(\Z)$
(this forces $k$ to be even),
and ``holomorphic at the cusp'' (see \S~\ref{subsec:prelim}). 
A form  is {\em cuspidal} if it vanishes at the cusp. 
For a modular form of weight $k$, let  $\nu(f)$
be the number of inequivalent zeros of $f$ in $\HH$,
with the convention that  a zero at $z$ is counted with weight $w(z)$ 
inverse to the number of elements of $SL_2(\Z)/\{\pm I\}$ fixing $z$. 
 Then $\nu(f)\leq k/12$. 

For a sequence of modular forms, where we assume that the number of
inequivalent zeros
tends to infinity, we would like to examine the manner in which the
resulting configuration of zeros is distributed in the modular domain
$SL_2(\Z)\backslash \HH$.

\subsection{Example: Eisenstein series}
% Take a fixed form $f_0$ and consider its powers
% $f_0^m$. These obviously have the same zeros as the original form, and
% so trivially $\ZF_\infty(f_0^m) = \ZF_\infty(f_0)$ and
% $\ZF(f^m)=\ZF(f)$.  As a special case take
% $f_0$ to be the modular discriminant, which is the unique cusp form of
% weight $12$, and  has has a simple zero at the
% cusp $i\infty$ and no other zeros. Then $f_0^m$ is a cusp form of
% weight $k=12m$ and $\ZF_\infty(f_0^m) =\frac{1}{k/12}m\delta_{i\infty} =
% \delta_{i\infty}$.
%
% A more interesting example are
The  Eisenstein series are non-cuspidal modular forms of weight $k>2$, 
given by the sum
$$
E_k(z)  = \frac 12\sum_{(c,d) = 1} \frac 1{(cz+d)^{k}}\;.
$$
%These are noncuspidal forms of weight $k$.
F.K.C. Rankin and  Swinnerton-Dyer \cite{Rankin-SD} showed that
all zeros (in the fundamental domain)
of $E_k$ lie on boundary arc $\{|z|=1\}$.
Moreover, as $k\to \infty$, the zeros become uniformly distributed 
on the unit arc. 
R.A. Rankin \cite{Rankin} gave a similar result for the
(cuspidal) Poincar\'e series
$$
G_k(z,m) = \sum_{\gamma\in \Gamma_\infty\backslash \Gamma} \frac
{\exp(2\pi i m\gamma z)}{(cz+d)^k}
$$
with $m\geq 1$ fixed, $k\gg m$. See \cite{Getz} for yet another
such example.

\subsection{Example: higher Weierstrass points}
An example in a slightly different setting are the higher Weierstrass
points on a (compact) Riemann surface $C$.
Let $S_k(C)$ be the space of holomorphic forms of weight $k$ on
$C$, $d:=\dim S_k(C)$, and let $f_1,\dots , f_d$  be a basis
for $S_k(C)$.
The Wronskian is defined as
$$
\Wr_k(\{f_j\};z):= \det ( f^{(i)}_j )_{1\leq i,j\leq d}
$$
It is a modular form of weight $d(k+d-1)$.

A {\em $k$-th order Weierstrass point} is a zero of $\Wr_k$.
%This definition is easily seen to be independent of the choice of basis
Equivalently, these are points of $C$ where there is a nonzero form
which vanishes to order $\geq \dim S_k(C)$.

It was conjectured by Bers, and proved by Olsen \cite{Olsen}
that as $k\to \infty$, these become dense in $C$. Mumford
(see page 11 of \cite{Mumford}) and Neeman \cite{Neeman} showed
that in fact these become equidistributed with respect to the Arakelov (or
Bergmann) measure on $C$, which arises from the metric on $C$ gotten
by embedding the curve  in its Jacobian and pulling back the flat
metric.

\subsection{}
In contrast to these examples,
we consider the case of cuspidal Hecke eigenforms. Recall that
the Hecke operators
$$T_n f(z):=\frac 1n\sum_{ad=n}a^k\sum_{b\mod d} f(\frac{az+b}d)$$
act on on the space of cusp forms of weight $k$,  
commute with each other and are self-adjoint
with respect to the Petersson inner product
$$
\langle f,g \rangle:=
\int_{SL_2(\Z)\backslash \HH} f(z)\overline{g(z)}y^k \frac{dxdy}{y^2}
$$
%on $S_k(SL_2(\Z))$.
Thus the space of cusp forms admits a basis consisting of joint
eigenfunctions of all Hecke operators (the Eisenstein series $E_k$ is
also an eigenfunction). The cuspidal Hecke eigenfunctions have a
simple zero at the cusp.
%so the bulk of their zeros lie in the finite part of the 
% fundamental domain.
For these Hecke eigenforms, we have:
\begin{thm}\label{thm:Hecke}
Assume the Generalized Riemann Hypothesis (GRH). Let $\{f_k\}$ be a
sequence of cuspidal Hecke eigenforms,
then as $k\to \infty$ their zeros are equidistributed in
$SL_2(\Z)\backslash \HH$
with respect to the normalized hyperbolic %Poincar\'e
measure  $dV(z) = \frac 3{\pi}\frac{dxdy}{y^2}$.
\end{thm}

%\begin{conj}\label{conjecture}%{\bf Conjecture:}
%{\em If $f_k\in S_k(SL_2(\Z))$ are cuspidal
%Hecke eigenforms,
%then as $k\to \infty$ their zeros are equidistributed in
%$SL_2(\Z)\backslash \HH$
%with respect to the normalized Poincar\'e measure
%$dV(z) = \frac 3{\pi}\frac{dxdy}{y^2}$.   }
%\end{conj}

%Recall that the Poincar\'e measure is associated to the constant
%negative curvature metric on the surface.
Equidistribution here means that for any nice {\em compact} subset
$\Omega \subset SL_2(\Z)\backslash \HH$,
the proportion of (inequivalent) zeros of such a Hecke eigenform
which lie in $\Omega$ is asymptotically the relative area
of $\Omega$:
$$
\frac 1{\nu(f)} \sum_{z\in \Omega} w(z)\ord_z(f)
\sim \int_\Omega dV(z)
$$

%As  key evidence for the validity of our conjecture, we have
%\begin{thm}
%  The Generalized Riemann Hypothesis implies
%  Conjecture~\ref{conjecture}.
%\end{thm}

\subsection{}
Theorem~\ref{thm:Hecke} is the outgrowth of some ideas from the 
theory of ``quantum chaos''. 
The mechanism is that equidistribution of the zeros of a sequence of
forms $f_k$  is implied by equidistribution of the ``masses''
$y^k|f_k(z)|^2 \dV(z)$ of the forms.
This idea was discovered by  Nonnenmacher and Voros \cite{NV}
in the context of quantum maps. 
Around the same time a very general result of this sort was obtained
by Shiffman and Zelditch \cite{SZ}
%by using the same argument,
for zeros of sections of high powers of a positive holomorphic
hermitian line bundle over any compact complex manifold
(\cite{NV} deal with curves of genus one).
%
%That is starting with a holomorphic line bundle over a complex
%projective manifold, which is assumed to be equipped with a hermitian
%metric $h$ which is positive
%
All these are in a compact setting.  
The analogous result in our case  is 
\begin{thm}\label{thm ergodicity}
Suppose that $\{f_k\}$ is a sequence of $L^2$-normalized
cusp forms ($f_k$ of weight $k$) 
for which the bulk of zeros lie in the fundamental domain:
$\nu(f_k) \sim k/12$.
Assume that for some $c>0$, we have
\begin{equation}\label{weak mass equid}
y^k|f_k(z)|^2 \dV(z)\toweak c\cdot \dV(z)
\end{equation}
%(where as before $dV(z) = dxdy/y^2$ is the Poincar\'e measure)
Then the zeros  of $f_k$ are equidistributed with respect to 
$\dV(z)$.
% $$
% \frac 1{\nu(f)}\sum_{j} w(z_{j,k}) \delta(z-z_{j,k})\toweak \dV(z)
% $$
\end{thm}
Here $\toweak$ denotes weak convergence when we test against compactly
supported functions.

%That there is such a relation is made plausible by the the following
%expression for the delta-masses of the zeros:
%$$
%\Delta\log|f_k(z)| = 2\pi \sum_{f_k(z_j)=0} \delta(z-z_j)
%$$
% In section~\ref{sec:argument} I will explain that a weakened form of
% the equidistribution of masses implies the equidistribution of the
% zeros:

 The proof of Theorem~\ref{thm ergodicity},  given in
Section~\ref{sec:argument}, follows closely  the argument of
\cite{NV,  SZ}, with care taken to handle the complications due to the
lack of compactness of the modular domain. 
We require the hypothesis  \eqref{weak mass equid}, which is weaker than 
equidistribution of the masses as it allows ``leakage'' of some of
 the mass at the cusp (which cannot happen in a compact setting). 
This slightly weaker version of equidistribution of masses is 
what  Lindenstrauss    \cite{Lindenstrauss} proved
 unconditionally  for the analogous  case of Maass  forms, though he
 cannot exclude $c=0$.   
In the holomorphic case we do not have an unconditional proof of 
\eqref{weak mass equid}. However it has been known for some time to 
follow (with $c=1$) from GRH \cite{AQC, T.Watson}.  
This is sketched in section~\ref{sketch of GRH}.  
%via a subconvexity bound for $L$-functions.  

%  We may   further weaken the condition by assuming weak convergence to
%  the absolutely continuous measure   $c(z)dV(z)$ where $c(z)\geq 0$
%  does not vanish on sets of positive measure and still keep the
%  conclusion.  In the case of Maass forms, this is a vacuous
%  generalization since any weak limit of the masses (viewed on the
%  unit cotangent bundle) has to be invariant under the geodesic flow,
%  as follows from Egorov's theorem, and so such $c(z)$ has to be
%  constant.  In the case of holomorphic forms I do not know of
%  an analogous statement.

%\subsection{Contents}
%Section~\ref{sec:PL} provides some background from potential
%theory, in particular giving a version of the well known relation
%between the delta function at the zeros of $f$ and the Laplacian of
%$\log|f|$. The proof of Theorem~\ref{thm ergodicity}  is given in
%Section~\ref{sec:argument}
%following closely  the argument of  \cite{NV,  SZ}, with care taken to
%handle the complications due to the lack of compactness of the modular
%domain.  
%Appendix~\ref{sec:bound} presents a bound for the supremum of
%$y^{k}|f(z)|^2$  as  $k\to\infty$.

{\bf Acknowledgements:} 
This manuscript is based on notes written in October 1999,
 while enjoying the   hospitality of the Institute of Advanced
 Studies in Princeton. I have benefited from several discussions with
 Peter Sarnak, Misha Sodin  and Steve Zelditch.  
%Partially supported by grant number 2002088  
%from the United States-Israel Binational Science
% Foundation (BSF), Jerusalem, Israel. 

%\input{prelim}
\section{Potential theory on $\Gamma\backslash \HH$}\label{sec:PL}

\subsection{Preliminaries} \label{subsec:prelim}
We review some definitions:
A modular form of weight $k$ for $SL_2(\Z)$
%for the modular group $SL_2(\Z)$
is a holomorphic function   on the upper half-plane $\HH$,
transforming as 
$$
f(\frac{az+b}{cz+d}) = (cz+d)^k f(z), \quad \mbox{for all }
\begin{pmatrix}  a&b\\c&d \end{pmatrix}\in SL_2(\Z)
$$
(this forces $k$ to be even),
and {\em holomorphic at the cusp}. % $i\infty$.
This means the following: Since $f(z)$ is
periodic, it can be expressed as a holomorphic
function $\tilde f(q)$ of $q=e^{2\pi iz}$ in the punctured disk
$0<|q|<1$. The requirement to be holomorphic at the cusp means that it
extends to a holomorphic function at $q=0$.
%This means that  as a function $\tilde f(q)$ of $q=e^{2\pi i z}$ it
%extends to a holomorphic function  at $q=0$.
A form  is {\em cuspidal} if it vanishes at the cusp, that is
$\tilde f(0)=0$.
The order of vanishing at the cusp $\ord_{\infty}(f)$
is defined as the order of vanishing at $q=0$ of $\tilde f(q)$.
%that is one expands $f$ in a power series in $q$: $f(z) =
%\sum_{n\geq 0} a_n q^n$   and then $\ord_{i\infty}(f)$ is the least
%$n\geq 0$ for which $a_n\neq 0$.

We denote by $\Gamma=SL_2(\Z)/\{\pm I\}$ and speak interchangeably
about modular forms for $\Gamma$ and for $SL_2(\Z)$.  
For each $z\in \HH$ let $\Gamma_z = \{\gamma\in \gamma: \gamma z=z \}$
be the stabilizer in $\Gamma$ of $z$, and set 
$$
w(z):=\frac 1{\#\Gamma_z} =
 \begin{cases}
 1/2,& z\in \Gamma i\\ 1/3,& z\in \Gamma e^{2\pi i/3}\\ 
1,& \mbox{ otherwise} \end{cases}
$$

For a modular form of weight $k$, let 
$$\nu(f)=\sum_{z\in SL_2(\Z)\backslash\HH } w(z)\ord_z(f)$$
be the (weighted) number of $\Gamma$-inequivalent zeros of $f$ in
$\HH$. Then $\nu(f)\leq \frac{k}{12}$, in fact  
$\nu(f)+ \ord_{\infty}(f) = \frac k{12}$. 
Note that $\nu(f)=0$ only for powers of the modular discriminant.

\subsection{}
% Let $\Gamma\subset PSL_2(\R)$ be a discrete subgroup of finite
% co-volume, and 
 For a smooth compactly supported function  $\phi\in C_c^\infty(\HH)$  on
the upper half-plane $\HH$  set 
$$ 
F_\phi(z):=\sum_{\gamma\in \Gamma} \phi(\gamma z) \in
C_c^\infty(\Gamma\backslash \HH)
$$
Note that 
$$\int_{\Gamma\backslash \HH} F_\phi(z) \frac{dxdy}{y^2} = \int_{\HH}
\phi(z)dxdy
$$ 
%For each $z\in \HH$ let $\Gamma_z = \{\gamma\in \gamma: \gamma z=z \}$
%be the stabilizer in $\Gamma$ of $z$, and set 
%$$w_\Gamma(z):=\frac 1{\#\Gamma_z}$$

Let $dV(z) = \frac 1{\vol(\Gamma\backslash \HH)}
\frac{dxdy}{y^2}$ be the normalized hyperbolic measure on the quotient
$\Gamma\backslash \HH$. 
\begin{lem}\label{dirac delta for zeros} 
Let $f$ be a (weakly holomorphic) modular form of weight $k$ for
$\Gamma$ and  let $\{z_j\}$ be a set of $\Gamma$-inequivalent zeros of
$f$ in $\HH$. Then 
\begin{multline}\label{eq:dirac delta for zeros} 
\sum_j w(z_j) F_\phi(z_j) = 
k\frac {\vol(\Gamma\backslash \HH)}{4\pi} 
\int_{\Gamma\backslash \HH} F_\phi(z)dV(z) \\
+  \frac 1{2\pi} \int_{ \HH} 
\log (y^{k/2}|f(z)|) \Delta \phi(z) dxdy 
\end{multline}
\end{lem}
\begin{proof}
Let $S=\{z\in \HH: f(z)=0\}$ be the set of zeros of the
form $f$.  It is a discrete set of points, which is
stable under $\Gamma$ ($\gamma S=S$ for all $\gamma\in \Gamma$). 
For any such $S$ we have
\begin{equation}\label{eq:app1} 
\sum_{z\in \Gamma\backslash S}  w(z) F_\phi(z) = \sum_{s\in S} \phi(s)
\end{equation}
Indeed, 
\begin{equation*}
  \begin{split}
\sum_{s\in S} \phi(s) &= 
\sum_{z\in \Gamma\backslash S} \sum_{s\in  \Gamma z}\phi(s) \\
& = \sum_{z\in \Gamma\backslash S} \sum_{\gamma\in \Gamma/\Gamma_z}
\phi(\gamma z) \\
&= \sum_{z\in \Gamma\backslash S} \frac  1{\#\Gamma_z} \sum_{\gamma\in
  \Gamma} \phi(\gamma z) \\
&= \sum_{z\in \Gamma\backslash S} \frac  1{\#\Gamma_z} F_\phi(z)
  \end{split}
\end{equation*}
as required. 

We recall that the electrostatic potential for a point charge in the
plane is $\frac 1{2\pi} \log |z|$, that is for 
$\phi\in C_c^\infty(\C)$ we have  
$$ 
\int_{\C} \frac 1{2\pi}\log|z| \Delta\phi(z) dxdy = \phi(0)
$$ 
Consequently, for the holomorphic function $f(z)$ on $\HH$ we
have 
\begin{equation}\label{eq:app2} 
\int_{\HH} \frac 1{2\pi}\log|f(z)| \Delta\phi(z) dxdy = 
\sum_{s: f(s)=0} \phi(s)
\end{equation} 
for all $\phi\in C_c^\infty(\HH)$ 
(with multiple zeros repeated). Thus 
\begin{multline*}
%\begin{split}
\sum_j w(z_j) F_\phi(z_j) = 
\frac 1{2\pi} \int_{\HH} \log |f(z)| \Delta \phi(z) dxdy\\
= - \frac 1{2\pi} \int_{\HH} \log y^{k/2} \Delta\phi(z)dxdy 
 + \frac 1{2\pi} \int_{\HH} \log(y^{k/2}|f(z)|) \Delta \phi(z) dxdy 
%\end{split} 
\end{multline*} 
The first term above  is transformed via integration by parts into 
$$
-\frac 1{2\pi}\int_{\HH} \Delta(\log y^{k/2}) \phi(z) dxdy  = 
\frac k{4\pi} \int_{\HH} \phi(z) \frac{dxdy}{y^2} 
$$
Note that 
$$
\int_{\HH} \phi(z) \frac{dxdy}{y^2} = \int_{\Gamma\backslash \HH}
F_{\phi}(z)\frac{dxdy}{y^2}  = 
\vol(\Gamma\backslash \HH) 
 \int_{\Gamma\backslash  \HH} F_\phi(z) dV(z)
$$
so that we get 
$$
-\frac 1{2\pi}\int_{\HH} \Delta(\log y^{k/2}) \phi(z) dxdy  = 
k \frac {\vol(\Gamma\backslash \HH)}{4\pi} 
\int_{\Gamma\backslash  \HH} F_\phi(z) dV(z)
$$ 
as required.
\end{proof}

Remark: 
We may reformulate this in $\Gamma$-invariant form as follows: 
Instead of the Euclidean Laplacian 
$\Delta= \frac{\partial^2}{\partial  x^2} + 
\frac{\partial^2}{\partial y^2}$ we use  the hyperbolic Laplacian
$$
\laplaceh =
y^2\left(\frac{\partial^2}{\partial x^2} + \frac{\partial^2}{\partial
  y^2} \right)
$$ 
which is $\Gamma$-invariant, and instead of $|f(z)|$ we use
$y^k|f(z)|^2$ which is $\Gamma$-invariant. Then
\begin{multline}\label{invariant form} 
\sum_j w(z_j) F_\phi(z_j) = k \frac {\vol(\Gamma\backslash \HH)}{4\pi} 
\int_{\Gamma\backslash \HH} F_\phi(z)dV(z) \\
+  \frac{\vol(\Gamma \backslash \HH)}{2\pi} 
\int_{\Gamma\backslash \HH} 
\log (y^{k/2}|f(z)|) \laplaceh F_\phi(z) \frac{dxdy}{y^2} 
\end{multline}
To derive \eqref{invariant form} from Lemma~\ref{dirac delta for zeros},  
we transform the second term in \eqref{eq:dirac delta for zeros} 
by noting that since $\frac{dxdy}{y^2}$
is $\Gamma$~-~invariant, we have
\begin{equation*}
\begin{split} 
 \frac 1{2\pi} \int_{\HH} \log(y^{k/2}|f(z)|) \laplaceh \phi(z)
\frac{dxdy}{y^2}  &= 
 \frac 1{2\pi} \int_{\Gamma\backslash \HH} \log(y^{k/2}|f(z)|)
 \sum_{\gamma\in \Gamma} (\laplaceh \phi)(\gamma z) \frac{dxdy}{y^2} \\
&= \frac 1{2\pi} \int_{\Gamma\backslash \HH} \log(y^{k/2}|f(z)|) 
F_{\laplaceh \phi}(z) \frac{dxdy}{y^2} \\
&= \frac{\vol(\Gamma \backslash \HH)}{2\pi} 
\int_{\Gamma\backslash \HH} \log(y^{k/2}|f(z)|) 
F_{\laplaceh \phi}(z) dV(z)
\end{split}
\end{equation*}
This proves \eqref{invariant form}  once we note that 
$F_{\laplaceh \phi}=\laplaceh F_\phi$ since $\laplaceh$ is
$\Gamma$-invariant.

\subsection{Background on sub-harmonic functions} \label{sec:sub-harmonic} 

%\subsection{}
Let $\Omega\subset \C=\R^2$ be an open connected set. 
An upper semi-continuous\footnote{$u:\Omega\to [-\infty,\infty)$ is
upper semi-continuous if  for all $c\in \R$, the set $\{z:u(z)<c\}$ is open.}
  function $u:\Omega\to [-\infty,\infty)$ is {\em sub-harmonic} if 
$u$ is not identically $-\infty$ and for all $z\in \Omega$,  
$$
u(z)\leq \frac 1{2\pi} \int_0^\infty u(z+re^{i\theta}) d\theta
$$
whenever the disc $\{w:|z-w|\leq r\}$ is contained in $\Omega$.

Note: $u$ is sub-harmonic  $\Leftrightarrow$ $\Delta u\geq 0$ as a
distribution. 

A fundamental example of a sub-harmonic function is $\log |f(z)|$ where
$f(z)$ is holomorphic in $\Omega$. 
Moreover we have 
$$
\Delta\log |f(z)| = 2\pi \sum_{f(z_j)=0} \delta(z-z_j)
$$
We denote the space of subharmonic functions on $\Omega$ by $\SHO$. 

%\subsection{}
A basic compactness property of sub-harmonic functions is
(c.f. \cite{Hormander}, Thm 4.1.9): 
\begin{lem}[Compactness property] \label{compactness thm}
If $\{u_j\}\subset \SHO $ are locally uniformly upper bounded (that is
for all compact $K\subset \Omega$ there is $c_K$ so that $u_j(z)\leq
c_K$ for all $z\in K$ and all $j$) then either 
\begin{enumerate} 
\item $u_j\to -\infty$ uniformly on compacta, OR 
\item There is a subsequence $u_{j_k}$ which converges weakly 
%, even in $L^1_{loc}(\Omega)$ 
to some $u\in \SHO$. 
\item \label{third}
Moreover, in this case $\limsup u_j \leq u$ and $\limsup u_j = u$
almost everywhere. 
\end{enumerate}
\end{lem}

%\begin{proof}
%Proof of part~\ref{third}: 
%Suppose we have a weakly convergent sequence $v_j\toweak v$. 

%\noindent i)  
%Firstly, the  limit $v$ is sub-harmonic (use $\Delta v\geq 0$). 

%\noindent ii) Moreover, $\limsup v_j\leq v$: 

%Indeed, taking an approximate
%identity\footnote{$\psi_\epsilon(x)=\frac 1{\epsilon^2}\psi(\frac
%x\epsilon)$, with $0\leq \psi\in C_c^\infty(\R^2)$, $\psi$ radial and
%$\int\psi(x)dx=1$.} $\psi_\epsilon$, 
%we have by weak convergence that  $\psi_\epsilon\conv v$ as $j\to \infty$. 
%Since $v_j$ are sub-harmonic, 
%$v_j\leq \psi_\epsilon\conv v_j\to \psi_\epsilon\conv v$. 
%Thus $v_j\leq \psi_\epsilon\conv v$ and taking $\epsilon\to 0$ gives
%$\limsup v_j\leq v$. 

%\noindent iii) We then have $\limsup v_j=v$ almost everywhere:  

%This we see by taking a positive test-function 
%$0\leq \chi\in C_c^\infty(\Omega)$  we have since $v_j\toweak v$ that 
%$\int \chi v=\lim_j \int \chi v_j$. By Fatou's lemma, $\lim_j\int \chi
%v_j \leq \int \chi \limsup v_j$, and by (ii), $\int\chi\limsup v_j
%\leq \int \chi v$. 
%Thus we find $\int\chi v\leq \int\chi \limsup v_j \leq \int\chi v$,
%and so we have equality $\int\chi\limsup v_j = \int \chi v$ for all
%$\chi$ as above. Therefore $\limsup v_j =v$ almost everywhere. 
%\end{proof}

%\subsection{} 
The following is known as ``Hartogs' Lemma'' (\cite{Ransford} Thm 3.4.3):
\begin{lem} \label{Hartogs' Lemma}
If $\{u_j\}\subset \SHO$ are locally uniformly bounded above, and
there is a continuous function $\phi:\Omega\to \R$ such that 
$\limsup u_j \leq \phi, then \max(u_j,\phi)\to\phi $ locally uniformly
on $\Omega$ as $j\to\infty$. 
\end{lem}

\section{Ergodicity of eigenfunctions implies equidistribution of zeros}
\label{sec:argument}

\subsection{Proof of Theorem~\ref{thm ergodicity}} 
Assume that there is some $c>0$ so that $y^k|f_k|^2 \dV\toweak c dV$
for a sequence of $k\to\infty$. 
We will show that  necessarily the zeros $\{z_j\}$ of $f_k$ become
equidistributed relative to $\dV$ as $k\to\infty$.  
We need to show that for all $F\in C_c^\infty(\Gamma\backslash \HH)$
we have
$$
\frac 1{\nu(f_k)} \sum_j w(z_j) F(z_j) \sim 
\int_{\Gamma\backslash\HH} F(z)dV(z), \qquad k\to \infty
$$
It suffices to show this for $F$ of the form $F(z)=F_\phi(z) =
\sum_{\gamma\in \Gamma}\phi(\gamma z)$, $\phi\in C_c^\infty(\HH)$, 
%\begin{verbatim}
%  Elaborate on this !
%\end{verbatim}
since these are dense in $C_c^\infty(\Gamma\backslash \HH)$ with
respect to the uniform topology.
\footnote{ To see this it suffices to show that for fixed
  ball $B\subset \HH$, the space $\{F_\phi:\supp \phi \subset B\}$ is
  dense in  $C(\pi(B))$, $\pi:\HH\to \Gamma\backslash \HH$ the
  projection, which  follows by applying the Stone-Weierstrass
  theorem. Note that the uniform closure of  the space
  $\{F_\phi:\phi\in C_c^\infty(\HH)\}$ is $C_0(\Gamma\backslash \HH)$, the
continuous functions vanishing at infinity.}  
For these we have the fundamental identity of 
Lemma~\ref{dirac delta for zeros} 
\begin{equation}
  \frac 1{\nu(f)} \sum_j w(z_j) F_\phi(z_j) =
\frac{k/12}{\nu(f_k)}   \int_{\Gamma\backslash\HH} F_\phi(z)dV(z) +
\mathcal E
\end{equation}
where 
$$
\mathcal E= 
\frac 1{2\pi \nu(f_k)} \int_{\HH} \log(y^k|f_k(z)|^2) \Delta\phi(z)dxdy 
$$
Since we assume that $\nu(f_k)\sim k/12$, the main term is the desired
one (the mean value of $F_\phi$) and we need to show that 
$\mathcal E\to 0$ as $k\to \infty$, i.~e.~ that for all $\phi \in
C_c^\infty(\HH)$, 
$$
\int_{\HH}  \frac 1k\log (y^k|f_k(z)|^2) \Delta\phi(z) dxdy \to 0
$$
or that 
\begin{equation}\label{false}
\int_{\HH} \frac 1k\log (|f_k(z)|^2)\Delta\phi(z) dxdy \to 
\int_{\HH} -\log y \Delta\phi(z) dxdy 
\end{equation} 

Assume that this is false, that is there is some test function
$ \phi_0\in C_c^\infty(\HH/\Gamma)$ and a subsequence of the $f_k$
(which we will continue to call $f_k$ for notational ease) for which
\eqref{false} fails. We will proceed to derive a contradiction. 

We begin by listing some properties of $v_k:=\frac 1k\log |f_k^2|$. 
These are sub-harmonic in $\HH/\Gamma$. By Proposition~\ref{bound prop}, we
have  $y^k|f_k(z)|^2 \ll k$ uniformly on compacta and so 
\begin{equation}
v_k\leq -\log y+\frac{\log k+O(1) }k
\end{equation}
thus
\begin{equation}\label{upper bound for limsup}
\limsup v_k \leq -\log y
\end{equation}
locally uniformly. 
Thus the family $\{v_k\}\subset \SH$ is locally bounded above. By the
compactness property (Lemma~\ref{compactness thm}), this gives us two
possibilities: 

i) $v_k\to-\infty$ locally uniformly, OR 

ii) $\{v_k\}$ has a weakly convergent subsequence. 

We will dispose of both possibilities. 

Option i) On the support of the test function $\phi_0$, we have
$v_k\to-\infty$ uniformly, so that there is some $K$ so that for all
$k\geq K$ and all $z\in \supp \phi_0$ we have 
$$
v_k(z)\leq -2H
$$
where $H=\max\{\im z:z\in \supp \phi_0 \}$. But then exponentiating
we find $$|f_k|^2\leq e^{-2kH}$$ so that for all $\phi$ supported
inside $\supp\phi_0$  (and so that $\supp \phi$ is contained in a
single fundamental domain)
$$\int_{\Gamma\backslash \HH} \phi(z)|f_k(z)|^2y^k \frac{dxdy}{y^2} \to 0$$
contradicting the
assumption that $y^k|f_k|^2 \dV\toweak c \dV$.

Option ii) We assume $\{v_k\}$ has a weakly convergent subsequence,
which for notational convenience we continue calling $\{v_k\}$, which
then converges to some $v\in \SH$, and moreover $\limsup v_k\leq v$  and
these are equal almost everywhere. Then by 
\eqref{upper bound for limsup}, $v+\log y\leq 0$ almost
everywhere. 

Since $v_k\toweak v$ but 
$\int_{\HH} v(z)\Delta\phi_0(z)dz \neq \int_{\HH} -\log y \Delta
\phi_0(z)dz$, we know that $v\neq -\log y$ on a set of positive
measure. Since $v$ is upper semi-continuous, there is some $\delta>0$
so that $v<-\log y-\delta$ on some nonempty open relatively compact
set $U$.  

Since $v=\limsup v_k<-\log y-\delta$, 
by Hartogs' lemma~\ref{Hartogs' Lemma}, there is $K=K(\delta,U)$ so
that for all $k\geq K$, $v_k\leq -\log y -\delta/2$ on $U$.
That is to say, 
$$
\frac 1k \log |f_k|^2<-\log y-\delta/2
$$
or 
$$
y^k|f_k|^2\leq e^{-k\delta/2}
$$
on $U$, which as in the first option contradicts 
$y^k|f_k|^2 \dV\toweak c \dV$.

\section{Sketch of the proof of Theorem~\ref{thm:Hecke} } 
\label{sketch of GRH}
As is well understood by now \cite{AQC, T.Watson}, the
equidistribution of masses \eqref{weak mass equid}  (with $c=1$!) 
follows from the Generalized Riemann Hypothesis
(for certain automorphic $L$-functions), in fact from a ``subconvexity
estimate'' for the central value of certain L-functions.  
To relate the equidistribution of the masses for Hecke eigenforms to GRH,  
one has to examine the behavior of the ``periods'' 
$$\int_{SL_2(\Z)\backslash\HH } g(z)|f_k(z)|^2 y^k dV(z)$$ 
where $g$ is  fixed and $k\to \infty$.  
%Hecke-Maass form, either an Eisenstein series or a cusp form.  
We need to show that whenever $g$ is orthogonal to the
constants then the period vanishes as $k\to \infty$. 

In the case where $g$ is  Eisenstein series 
%was  treated unconditionally by Luo and Sarnak \cite{LS}. In that
%case 
the period  is the (completed) Rankin-Selberg $L$-function 
$L^*(1/2+it,f_k\times f_k)$, the standard convexity estimates give a
bound of $k^\epsilon$ (for all $\epsilon>0$)   
and vanishing as $k\to \infty$ follows from GRH (see \cite[\S 4]{AQC}).

To handle the case when $g$ is a cuspidal Hecke-Maass form,  
one exploits  Watson's formula \cite{T.Watson} which relates 
the period with the central value of the triple product $L$-function:  
$$
L(\frac 12,g\times f_k\times f_k) 
$$
%\begin{verbatim} divide by sym^2 \end{verbatim}
%for $c>0$ ($f_k$ and $g$ are $L^2$-normalized). 
%and where $L^*$ denotes the completed $L$-function. 
To get the decay as $k\to \infty$, 
it transpires that, again, one has to beat the ``convexity bound'' on the
central value, for which one uses GRH. 

%The equivalence of non-vanishing for the central value of the triple
%product $L$-function and that of the periods was conjectured by
%Jacquet and proven by Harris and Kudla \cite{Kudla-Harris}. 

In the case of CM-forms\footnote{These only exist for certain
congruence subgroups of the modular group},  Sarnak \cite{SarnakJFA}  
proved %(without going via $L$-functions) 
that their masses are equidistributed. 
Thus the reasoning above implies (unconditionally) that the zeros of 
CM-forms are equidistributed w.r.t. $\dV$ as $k\to\infty$.

%\input{refs}

%\newpage 

\appendix 

\section{An $L^\infty$-bound for cusp forms}\label{sec:bound} 

A holomorphic form of weight $k$ for $SL_2(\Z)$ is bounded, since we
require that the $q$-expansion contains no negative powers of $q$. 
We will need a bound 
on the supremum with an explicit dependence on the weight $k$. 
For our purposes any bound on $y^{k/2}|f(z)|$ 
which is sub-exponential  in $k$ for fixed $z$ suffices. 
Below we derive a bound of size $k^{1/2}$; 
since the dimension of the space of cusp forms grows linearly in $k$,  
this bound is optimal (for arbitrary forms) as far as the $k$-dependence. 
The proof is adapted from \cite[Lemma A.1]{IS}  who treat the case of
Maass forms, though unlike \cite{IS} we do not need to control the
$z$-aspect.  

\subsection{The incomplete Gamma function} 
We first recall some properties of the  incomplete Gamma function
$\Gamma(a,x)$ defined for $x>0$ by 
$$
\Gamma(a,x) := \int_x^\infty e^{-t}t^{a}\frac {dt}t
$$
  From the definition it is clear that $\Gamma(a,x)$ 
is decreasing in $x$. 
For integer $k$ we have 
\begin{equation}\label{evaluation of incg}
\Gamma(k-1,x) = %\int_x^\infty e^{-t}t^{k-1}\frac {dt}t =
(k-2)!e^{-x}\sum_{m=0}^{k-2} \frac{x^m}{m!}
\end{equation}

What is not completely straightforward is the asymptotic behavior of
$\Gamma(a,x)$ when $a,x$ tend to infinity. 
We will need the asymptotic: 
\begin{equation}\label{ramanujan}
\Gamma(k-1,k) \sim \frac 12\Gamma(k-1),\qquad k\to \infty
\end{equation}
%\vskip .5cm
%\noindent{\bf Remark:} On the relation
%$\Gamma(k-1,k) \sim \frac 12\Gamma(k-1)$: 
By \eqref{evaluation of incg}, 
% Since $\Gamma(k-1,k) = (k-2)!e^{-k}\sum_{m=0}^{k-2}\frac{k^m}{m!}$,
the relation \eqref{ramanujan} is equivalent to
\begin{equation}\label{weak 1}
e^{-k}\sum_{m=0}^{k-2}\frac{k^m}{m!}\sim \frac 12
\end{equation}
This is close to a conjecture of Ramanujan
\cite{Ram1, Ram2, Szego, G.N.Watson} that
$$
e^{-k}\left( \sum_{m=0}^{k-1}\frac{k^m}{m!} +
\theta\frac{k^k}{k!}\right) 
= \frac12
$$
for some $\frac 13\leq \theta=\theta(k)\leq \frac  12$;
%presumably 
he showed
\cite{Ram2} that $\theta(\infty)= 1/3$.
We need  the weaker asymptotic \eqref{weak 1} or equivalently (by
Stirling's formula), that %$\theta(k) = o(\sqrt{k})$, or that
\begin{equation}\label{weak 2}
e^{-k}\sum_{m=0}^{k}\frac{k^m}{m!} \sim \frac 12
\end{equation}
Here is a quick ``proof'' using the central limit theorem (this kind
of argument goes back to Mark Kac \cite{Kac}):
Let $X_1,\dots, X_k$ be independent Poisson variables with parameter
$1$ (so having mean and variance $1$).
Let $S_k=X_1+\dots +X_k$ be their sum. By the
Central Limit Theorem,    $\frac{S_k-k}{\sqrt{k}}$ is asymptotically
normal, in particular as $k\to \infty$,
$$
\mbox{Prob}(\frac{S_k-k}{\sqrt{k}}\leq 0) \to \frac 12
$$
that is
$$
\mbox{Prob}(S_k\leq k) \to \frac 12
$$
However since $S_k$ is Poisson with parameter $k$, we have
$$
\mbox{Prob}(S_k\leq k) = e^{-k}\sum_{m=0}^k\frac{k^m}{m!}
$$
which gives \eqref{weak 2}.

\subsection{The supremum of cusp forms}
\begin{prop}\label{bound prop}
Let $f$ be a cusp form of weight $k$ for $SL_2(\Z)$.
%Then for $z$ in the standard fundamental domain, and $k\to \infty$,
%$$
%
%$$
%
%That is
Then uniformly for $z$ in compact subsets of $\HH$, 
$$
\frac {y^k|f(z)|^2}{\langle f,f\rangle} \ll  k 
$$
where $\langle f,f \rangle$ is the Peterson inner product. 
\end{prop}

\begin{proof} 
Consider the integral of $y^k|f(z)|^2$ over the ``Siegel set'' 
$$\mathfrak S:=\{z=x+iy:0\leq x\leq 1,\quad y>1/4\pi \}$$ 
This integral converges
since $f$ is exponentially decreasing in the cusp. 
Since the Siegel set $\mathfrak S$  %$\{z: \Im(z)\geq Y \}$ 
is contained in a fixed number of translates 
of the standard  fundamental domain 
$\mathcal F=\{z=x+iy: |z|\geq 1, |x|\leq 1/2\}$ 
we have
\begin{equation*}
\int_{\mathfrak S} |f(z)|^2 y^k\frac{dxdy}{y^2} \ll
\langle f, f\rangle
\end{equation*}

On the other hand, using the Fourier expansion $f(z) = \sum_{n\geq 1}
a(n)e^{2\pi i nz}$ and Parseval, we find
\begin{equation*}
\begin{split}
\int_{\mathfrak S} |f(z)|^2 y^k\frac{dxdy}{y^2}  &=
\sum_{n\geq 1} |a(n)|^2 \int_{1/4\pi}^\infty e^{-4\pi ny} y^{k-2}dy \\
&=
\sum_{n\geq 1}\frac{|a(n)|^2 }{(4\pi n)^{k-1}} \Gamma(k-1, n)
\end{split}
\end{equation*}
Thus we find that
\begin{equation}\label{2 norms}
\sum_{n\geq 1}\frac{|a(n)|^2 }{(4\pi n)^{k-1}} \Gamma(k-1, n) 
\ll \langle f, f\rangle
\end{equation}

%Note: Since $k-1$ is an integer, we have
%$$
%\Gamma(k-1,x) = (k-2)!e^{-x}\sum_{m=0}^{k-2} \frac{x^m}{m!}
%$$
%and therefore ommiting all terms in the sum except that with $m=k-2$
%gives
%$$
%\Gamma(k-1,x)\geq e^{-x}x^{k-2}
%$$

To estimate $y^k|f(z)|^2$, use the Fourier expansion and
Cauchy-Schwartz:  %with weights
\begin{equation*}
\begin{split}
y^k|f(z)|^2 &= y^k\left |\sum_{n\geq 1} a(n) e^{-2\pi ny} e(nx)
\right|^2 \\
&\leq  \sum_{n\geq 1}\frac{|a(n)|^2 }{(4\pi n)^{k-1}}
\Gamma(k-1, n) \times y^k \sum_{n\geq 1} e^{-4\pi n y}
\frac{(4\pi n )^{k-1}}{\Gamma(k-1,n)}
\end{split}
\end{equation*}

For the first sum, use \eqref{2 norms} to bound it by
$\langle f,f\rangle$.
Thus we find
$$
\frac{y^k|f(z)|^2}{\langle f,f\rangle} \ll
y^k \sum_{n\geq 1} e^{-4\pi n y}
\frac{(4\pi n )^{k-1}}{\Gamma(k-1, n)}
$$

To bound  the  sum, we split it up into ``small'' and ``large''
$n$'s.
%$$
%\sum_{n\geq 1} e^{-4\pi n y}
%\frac{(4\pi n y)^{k-1}}{\Gamma(k-1, 4\pi nY)} = I+II
%$$
For ``small'' $n$, that is those  satisfying $n\leq k$,
we can give an upper bound  as follows: First use
$
\Gamma(k-1,n)\geq \Gamma(k-1,k)$, for $n\leq k$.  
This gives
\begin{equation*}
\begin{split}
%I& =
 y^k \sum_{n\leq k} e^{-4\pi n y}
\frac{(4\pi n )^{k-1}}{\Gamma(k-1,n)}
& \leq \frac {y^k}{\Gamma(k-1,k) }
\sum_{n\leq k}  e^{-4\pi n y} (4\pi n )^{k-1} \\
&\leq
\frac {y}{\Gamma(k-1,k) } \sum_{n=1}^\infty
 e^{-4\pi n y} (4\pi n y)^{k-1}
\end{split}
\end{equation*}
We have\footnote{If $f(t)$  is increasing for $t<t_0$ and decreasing
for $t>t_0$ then $\sum_{n=1}^\infty f(n)\leq \int_0^\infty f(t)dt
+2f(t_0)$}
\begin{equation*}
\begin{split}
 \sum_{n=1}^\infty  e^{-4\pi n y} (4\pi n y)^{k-1}
&\leq
\int_0^\infty  e^{-4\pi x y} (4\pi x y)^{k-1} dx  +2e^{-(k-1)}
(k-1)^{k-1} \\
&\ll  \Gamma(k)(\frac 1y+ \frac 1{\sqrt{k}}) %{4\pi y} (1+O(k^{-1/2}) )
\end{split}
\end{equation*}
Thus  we find 
$$
y^k\sum_{n\leq k}   e^{-4\pi n y}
\frac{(4\pi n )^{k-1}}{\Gamma(k-1, n)}
\ll \frac {y \Gamma(k)} {\Gamma(k-1,k)}   (\frac{1}{y}+\frac 1{\sqrt{k}}) 
%{ \frac 12\Gamma(k-1)} 
%%\approx k %%\frac ky
$$
which by \eqref{ramanujan} is 
$$\ll k( 1+\frac y{\sqrt{k}} )$$

%Note: I believe this sum should be bounded by $k$ rather than
%$k^{3/2}$ !!

 For the sum over the remaining $n$'s, use \eqref{evaluation of incg} 
to get 
$$
\Gamma(k-1,n)\geq e^{-n} n^{k-2}
$$
which gives
\begin{equation*}
\begin{split}
 y^k \sum_{n> k} e^{-4\pi n y}
\frac{(4\pi n )^{k-1}}{\Gamma(k-1,n)}
 &\leq y^k \sum_{n>k} e^{-4\pi n y}
\frac{(4\pi n )^{k-1}}{e^{-n}n^{k-2} } \\
&\leq (4\pi y)^{k}  \sum _{n>k} e^{-n(4\pi y-1)} n\\
&\ll 
(4\pi y)^{k} \frac {(k+1) e^{-(k+1)(4\pi y-1)} }{(1-e^{-(4\pi
    y-1)})^2} \\
& \ll  k e^{-4\pi y} (\frac{4\pi ye}{e^{4\pi y}} )^k
%\frac{y^{k-1}}{y-Y} ke^{-k(\frac yY-1)}
\end{split}
\end{equation*}
which, for fixed $y\geq \sqrt{3}/2$, is negligible.
\end{proof}

%\begin{verbatim}
%Actually need to worry about the possibility $y=Y=\sqrt{3}/2$,
%that is the bottom elliptic points !
%\end{verbatim}

\end{document}